\documentclass{amsart}%
\usepackage{amsfonts}%
\usepackage{amsmath}%
\usepackage{amssymb}%
\usepackage{graphicx}
\newtheorem{theorem}{Theorem}
\theoremstyle{plain}

\newtheorem{proposition}{Proposition}

\numberwithin{equation}{section}
\begin{document}
\title{On a conjecture of Schroeder and Strake}

\begin{abstract}
We prove some rigidity results for compact manifolds with boundary. In
particular for a compact Riemannian manifold with nonnegative Ricci curvature
and simply connected mean convex boundary, it is shown that if the sectional
curvature vanishes on the boundary, then the metric must be flat.

\end{abstract}
\subjclass{53C24, 53C21}
\keywords{rigidity, nonnegative Ricci curvature, mean convex boundary, Reilly's formula}
\author{Fengbo Hang}
\address{Department of Mathematics, Princeton University, Fine Hall, Washington Road,
Princeton, NJ\ 08544}
\email{fhang@math.princeton.edu}
\author{Xiaodong Wang}
\address{Department of Mathematics, Michigan State University, East Lansing, MI 48824}
\email{xwang@math.msu.edu}
\thanks{F. Hang is partially supported by NSF grant DMS-0501050 and a Sloan Research Fellowship.}
\thanks{X. Wang is partially supported by NSF grant DMS-0505645.}
\maketitle

In \cite[Theorem 1]{SS}, Schroeder and Strake proved the following rigidity theorem.

\begin{quotation}
Let $\left(  M,g\right)  $ be a compact Riemannian manifold with convex
boundary and nonnegative Ricci curvature. Assume that the sectional curvature
is identically zero in some neighborhood $U$ of $\partial M$ and that one of
the following conditions holds:
\end{quotation}

\begin{itemize}
\item $\partial M$ is simply connected

\item $\dim\partial M$ is even and $\partial M$ is strictly convex at some
point $p\in\partial M$.
\end{itemize}

\begin{quotation}
Then $M$ is flat.
\end{quotation}

As remarked in \cite{SS}, the condition that the metric is flat in a whole
neighborhood of $\partial M$ is very strong. They conjectured that it suffices
to only assume that the sectional curvature vanishes on $\partial M$ and
proved this in the special case of a convex metric ball. The problem was
studied by Xia in \cite{X1, X2} who confirmed the conjecture under various
additional conditions: like the boundary has constant mean curvature or
constant scalar curvature, or the second fundamental form satisfies some
pinching condition etc. We refer to \cite{X1,X2} for the precise statements.
Here we present some results related to the conjecture.

\begin{theorem}
\label{thmmain} Let $M$ be a smooth compact connected Riemannian manifold with
boundary and nonnegative Ricci curvature. If every component of $\partial M$
is simply connected and has nonnegative mean curvature and the sectional
curvature of $M$ vanishes on $\partial M$, then $M$ is flat and $\partial M$
has only one component.
\end{theorem}

Therefore when $\partial M$ is simply connected the conjecture of Schroeder
and Strake is true. Moreover one only needs $\partial M$ to be mean convex
instead of convex. We remark that the conclusion that $\partial M$ has only
one component follows from theorems in \cite{I,K}. Below we will present a
different argument for it based on the Reilly's formula (\cite{Re}).

To continue the discussion we need to fix some notations. We will often write
$\left\langle \ ,\ \right\rangle $ for the metric on $M$ and denote the
connection as $D$. For convenience we write $\Sigma=\partial M$ and denote the
Levi-Civita connection and curvature tensor etc. of the induced metric on
$\Sigma$ as standard notations with a subscript $\Sigma$. Let $\nu$ be the
unit outer normal vector. The shape operator is given by $A\left(  X\right)
=D_{X}\nu$ and the second fundamental form is given by $h\left(  X,Y\right)
=\left\langle A\left(  X\right)  ,Y\right\rangle =\left\langle D_{X}%
\nu,Y\right\rangle $, here $X,Y\in T\Sigma$. The mean curvature
$H=\operatorname{tr}A$. Recall Reilly's formula (\cite[formula (14)]{Re}) for
a smooth function $u$ on $M$%
\begin{align*}
\frac{1}{2}\int_{M}\left(  \left(  \Delta u\right)  ^{2}-\left\vert
D^{2}u\right\vert ^{2}\right)  d\mu & =\frac{1}{2}\int_{M}Rc\left(  \nabla
u,\nabla u\right)  d\mu+\int_{\Sigma}\Delta_{\Sigma}u\cdot\frac{\partial
u}{\partial\nu}dS\\
& +\frac{1}{2}\int_{\Sigma}H\left(  \frac{\partial u}{\partial\nu}\right)
^{2}dS+\frac{1}{2}\int_{\Sigma}\left\langle A\left(  \nabla_{\Sigma}u\right)
,\nabla_{\Sigma}u\right\rangle dS.
\end{align*}

A special case of theorems in \cite{I,K} claims that if $M^{n}$ is a compact
connected Riemannian manifold with mean convex boundary $\Sigma$ and
nonnegative Ricci curvature, then $\Sigma$ has at most two components;
moreover if $\Sigma$ has two components, then $M$ is isometric to
$\Gamma\times\left[  0,a\right]  $ for some connected compact Riemannian
manifold $\Gamma$ with nonnegative Ricci curvature and $a>0$. For Theorem
\ref{thmmain}, it is clear $M$ can not have the product metric, hence $\Sigma$
has one component. It is interesting that one may give an argument for the
above special case based on Reilly's formula. Indeed, assume $\Sigma$ is not
connected, fix a component $\Sigma_{0}$ of $\Sigma$, then we may solve the
Dirichlet problem%
\[
\left\{
\begin{array}
[c]{l}%
\Delta u=0\text{ on }M,\\
\left.  u\right\vert _{\Sigma_{0}}=0,\\
\left.  u\right\vert _{\Sigma\backslash\Sigma_{0}}=1.
\end{array}
\right.
\]
Applying the Reilly's formula to $u$, we get%
\[
-\int_{M}\left\vert D^{2}u\right\vert ^{2}d\mu=\int_{M}Rc\left(  \nabla
u,\nabla u\right)  d\mu+\int_{\Sigma}H\left(  \frac{\partial u}{\partial\nu
}\right)  ^{2}dS.
\]
Hence $D^{2}u=0$. This implies $\left\vert \nabla u\right\vert \equiv c>0$.
Since $\nabla u=-c\nu$ on $\Sigma_{0}$ and $\nabla u=c\nu$ on $\Sigma
\backslash\Sigma_{0}$, we see $D_{X}\nu=0$ for $X\in T\Sigma$ i.e. $\Sigma$ is
totally geodesic. If we look at the flow generated by $\frac{\nabla u}{c}$,
then it sends $\Sigma_{0}$ to $\Sigma\backslash\Sigma_{0}$ at time $\frac
{1}{c}$ and hence $\Sigma$ has exactly two components. Note that the flow
lines are just geodesics. If we fix a coordinate on $\Sigma_{0}$, namely
$\theta^{1},\cdots,\theta^{n-1}$, let $r=\frac{u}{c}$, then we have
$g=dr\otimes dr+g_{ij}\left(  r,\theta\right)  d\theta^{i}\otimes d\theta^{j}%
$. Using $D^{2}r=0$, we see $\partial_{r}g_{ij}\left(  r,\theta\right)  =0$.
Hence $M$ is isometric to $\Sigma_{0}\times\left[  0,\frac{1}{c}\right]  $.

Under the assumption of Theorem \ref{thmmain} that the sectional curvature of
$M$ vanishes on $\Sigma$, it follows from Gauss and Codazzi equations that%
\begin{align*}
R_{\Sigma}\left(  X,Y,Z,W\right)   & =h\left(  X,Z\right)  h\left(
Y,W\right)  -h\left(  X,W\right)  h\left(  Y,Z\right)  ,\\
\left(  D_{\Sigma}\right)  _{X}h\left(  Y,Z\right)   & =\left(  D_{\Sigma
}\right)  _{Y}h\left(  X,Z\right)  ,
\end{align*}
where $X,Y,Z$ and $W$ belong to $T\Sigma$. By the fundamental theorem for
hypersurfaces \cite[part (2) of Theorem 21 on p63]{Sp} and the fact $\Sigma$
is simply connected, we may find a smooth isometric immersion $\phi
:\Sigma\rightarrow\mathbb{R}^{n}$ such that the second fundamental form of the
immersion $h_{\phi}=h$. If $\Sigma$ is convex then $\phi$ is an embedding by a
Hadamard type theorem of Sacksteder \cite{S}. With this immersion $\phi$ at
hands, Theorem \ref{thmmain} follows from the following proposition.

\begin{proposition}
Assume $M^{n}$ is a smooth compact connected Riemannian manifold with
connected boundary $\Sigma=\partial M$ and $Rc\geq0$. If $\phi:\Sigma
\rightarrow\mathbb{R}^{l}$ is an isometric immersion with $\left\vert H_{\phi
}\right\vert \leq H$ on $\Sigma$, here $H_{\phi}$ is the mean curvature vector
of the immersion $\phi$, then $M$ is flat. If moreover $\phi$ is an imbedding,
then $M$ is isometric to a domain in $\mathbb{R}^{n}$.
\end{proposition}

This is a generalization of a theorem of Ros \cite[Theorem 2]{R}, who derived
a congruence theorem for hypersurface in Euclidean space. Following the
argument of Ros, we will show by Reilly's formula that the harmonic extension
of the map $\phi$ is in fact an isometric immersion.

\begin{proof}
We may find a smooth function $F:M\rightarrow\mathbb{R}^{l}$ such that%
\[
\left\{
\begin{array}
[c]{c}%
\Delta F=0\text{ in }M;\\
\left.  F\right\vert _{\Sigma}=\phi.
\end{array}
\right.
\]
Applying the Reilly's formula to each component of $F$ and sum up we get%
\begin{align*}
& -\frac{1}{2}\sum_{\alpha}\int_{M}\left\vert D^{2}F^{\alpha}\right\vert
^{2}d\mu\\
& =\frac{1}{2}\int_{M}\sum_{\alpha}Rc\left(  \nabla F^{\alpha},\nabla
F^{\alpha}\right)  d\mu+\int_{\Sigma}\sum_{\alpha}\Delta_{\Sigma}\phi^{\alpha
}\cdot\frac{\partial F^{\alpha}}{\partial\nu}dS\\
& +\frac{1}{2}\int_{\Sigma}\sum_{\alpha}H\left(  \frac{\partial F^{\alpha}%
}{\partial\nu}\right)  ^{2}dS+\frac{1}{2}\int_{\Sigma}\sum_{\alpha
}\left\langle A\left(  \nabla_{\Sigma}\phi^{\alpha}\right)  ,\nabla_{\Sigma
}\phi^{\alpha}\right\rangle dS.
\end{align*}
Note%
\begin{align*}
\sum_{\alpha}\left\langle A\left(  \nabla_{\Sigma}\phi^{\alpha}\right)
,\nabla_{\Sigma}\phi^{\alpha}\right\rangle  & =\left\langle Ae_{i}%
,e_{j}\right\rangle e_{i}\phi^{\alpha}\cdot e_{j}\phi^{\alpha}=\left\langle
Ae_{i},e_{j}\right\rangle \phi_{\ast}e_{i}\cdot\phi_{\ast}e_{j}\\
& =\left\langle Ae_{i},e_{j}\right\rangle \delta_{ij}=\operatorname*{tr}A=H,
\end{align*}
here $e_{1},\cdots,e_{n-1}$ is a local orthonormal frame on $\Sigma$, hence%
\begin{align*}
0  & =\frac{1}{2}\sum_{\alpha}\int_{M}\left\vert D^{2}F^{\alpha}\right\vert
^{2}d\mu+\frac{1}{2}\int_{M}\sum_{\alpha}Rc\left(  \nabla F^{\alpha},\nabla
F^{\alpha}\right)  d\mu+\int_{\Sigma}H_{\phi}\cdot F_{\ast}\nu dS\\
& +\frac{1}{2}\int_{\Sigma}H\left\vert F_{\ast}\nu\right\vert ^{2}dS+\frac
{1}{2}\int_{\Sigma}HdS\\
& \geq\frac{1}{2}\sum_{\alpha}\int_{M}\left\vert D^{2}F^{\alpha}\right\vert
^{2}d\mu+\frac{1}{2}\int_{\Sigma}H\left(  \left\vert F_{\ast}\nu\right\vert
^{2}-2\left\vert F_{\ast}\nu\right\vert +1\right)  dS.
\end{align*}
Hence $D^{2}F^{\alpha}=0$ for all $\alpha$. It follows that $F^{\ast
}g_{\mathbb{R}^{l}}$ is parallel on $M$. We may find some $p\in\Sigma$ such
that $\left\vert H_{\phi}\right\vert >0$ at $p$, hence $H\left(  p\right)
>0$. From the argument above this implies $\left\vert F_{\ast}\nu\right\vert
=1$ at $p$ and $F_{\ast}\nu$ is perpendicular to $\phi_{\ast}\Sigma_{p}$,
hence $F^{\ast}g_{\mathbb{R}^{l}}=g_{M}$ at $p$. It follows that $F^{\ast
}g_{\mathbb{R}^{l}}=g_{M}$ on $M$, that is, $F$ is an isometric immersion and
$M$ is flat. Now assume $\phi$ is an imbedding. Let $\overline{D}$ be the
connection on $\mathbb{R}^{l}$, then $\overline{D}_{X}F_{\ast}Y-F_{\ast}%
D_{X}Y=XYF-\left(  D_{X}Y\right)  F=0$, it follows that $F:M\rightarrow
\mathbb{R}^{l}$ is a totally geodesic submanifold, hence the image lies in a
$n$ dimensional affine subspace. Without losing of generality we may assume
$l=n$ and $\Sigma$ is a compact hypersurface in $\mathbb{R}^{n}$, then there
exists a bounded open domain $\Omega$ such that $\partial\Omega=\Sigma$. Since
$F$ is an immersion, we see $F\left(  M\right)  \backslash\overline{\Omega}$
is both open and closed in $\mathbb{R}^{n}\backslash\overline{\Omega}$, hence
it must be empty. Based on this we may show $F:M\rightarrow\overline{\Omega}$
is a covering map and hence it must be a diffeomorphism.
\end{proof}

If we assume that $\partial M$ is convex, then it is clear from the above
discussion that $M$ is isometric to a convex domain in $%
\mathbb{R}
^{n}$. In fact in this case one may replace the nonnegativity of the Ricci
curvature by the much weaker nonnegativity of the scalar curvature, at least
when $M$ is spin.

\begin{theorem}
Let $M$ be a smooth compact connected Riemannian manifold with boundary and
nonnegative scalar curvature. If $M$ is spin, each component of $\partial M$
is convex and simply connected and the sectional curvature of $M$ vanishes on
$\partial M$, then $M$ is isometric to a convex domain in $%
\mathbb{R}
^{n}$.
\end{theorem}

\begin{proof}
For every component $\Gamma$ of $\partial M$, we have an isometric embedding
$\phi:\Gamma\rightarrow%
\mathbb{R}
^{n}$ which has $h$ as the second fundamental form. Let $\Omega$ be the convex
domain enclosed by $\phi\left(  \Gamma\right)  $. We glue $M$ and $%
\mathbb{R}
^{n}\backslash\Omega$ along $\Gamma$ via the diffeomorphism $\phi$ for all the
$\Gamma$'s and obtain a complete Riemannian manifold $N$ which has nonnegative
scalar curvature and is flat outside a compact set. Notice that the metric is
$C^{1}$ along the gluing hypersurface. Since $M$ is spin, we conclude by the
generalized positive mass theorem proved in \cite[Theorem 3.1]{ST} that $N$ is
isometric to $%
\mathbb{R}
^{n}$. It follows that $M$ is isometric to a convex domain in $%
\mathbb{R}
^{n}$.
\end{proof}

\end{document}